\documentclass[11pt]{article}
\usepackage{amsmath,amssymb,amsthm,eucal, graphicx, units}

\usepackage[utf8]{inputenc} 
\usepackage[T1]{fontenc}    

\usepackage{url}            
\usepackage{booktabs}       
\usepackage{amsfonts}       
\usepackage{nicefrac}       
\usepackage{microtype}      
\usepackage{lipsum}		
\usepackage{enumitem}
\usepackage{graphicx}

\usepackage[square,numbers]{natbib}
\usepackage{doi}
\usepackage{hyperref}
\usepackage{float}
\newtheorem{theorem}{Theorem}[section]
\newtheorem{corollary}{Corollary}[theorem]

\usepackage{subcaption}
\usepackage{pdfpages}
\usepackage{framed}
\usepackage{pdfpages}
\usepackage{color}

\usepackage[headings]{fullpage}
\usepackage{tocloft}
\usepackage{soul}
\usepackage{yfonts}
\usepackage[all]{xy}
\usepackage{float}
\usepackage{graphicx}
\usepackage{caption}
\usepackage{tkz-fct}
\usepackage{tikz}
\usetikzlibrary{braids}
\usetikzlibrary{calc,through,intersections, arrows, decorations.markings, matrix,shapes,snakes}

\newtheorem{definition}[theorem]{Definition}

\newtheorem{example}[theorem]{Example}

\newtheorem{question}[theorem]{Question}

\newtheorem{definition/proposition}[theorem]{Definition/Proposition}

\usepackage{cleveref}

\title{ Bridge Positions and Plat Presentations of Links}
\date{}
\author{ Seth Hovland}

\hypersetup{
pdftitle={Bridge Presentations and Plats},
}

\begin{document}
\maketitle

\begin{abstract}
This paper investigates the relationship between links in bridge position and plat presentations.  
We prove that Hilden double coset classes of plat presentations correspond exactly to bridge positions of a link up to bridge isotopy.  This correspondence allows us to translate algebraic questions about plat presentations into geometric questions about bridge positions.  
Using this perspective, we re-establish several results from  both frameworks.  In particular, we show that there is a unique Hilden double coset class for the $n$-bridge unknot in $S^3$, and that torus knots admit a single double coset class in plat position.  Finally, we discuss how this correspondence can be used to study plat closures of braids, which is the subject of ongoing research.

\end{abstract}


\section{Introduction}\label{sec:Introduction}

Bridge presentations and plat presentations are two classical ways of encoding links in $S^3$. 
Bridge presentations arise from Morse-theoretic considerations with respect to a height function, while plat presentations arise algebraically from braid closures with a fixed pairing of endpoints. 

These constructions are well known to be related: every bridge presentation can be isotoped into a plat form, and conversely every plat defines a bridge presentation.  The precise relationship between their natural equivalence relations is less explicitly recorded in the literature.

The main purpose of this paper is to identify the correct equivalence relation on plat presentations that corresponds exactly to bridge isotopy.

\begin{theorem}[Main Theorem]\label{maintheorem_intro}
The set of bridge isotopy classes of $n$-bridge presentations of links in $S^3$ is in bijection with Hilden double cosets in the braid group $B_{2n}$. 
Equivalently, two plat presentations represent bridge isotopic links if and only if their corresponding braids lie in the same Hilden double coset.
\end{theorem}

This result provides a direct correspondence between the geometric notion of bridge isotopy and the algebraic notion arising from braid groups. In particular, it allows questions about bridge presentations to be translated into questions about double cosets in the braid group with respect to the Hilden subgroup. Indeed, recent work by Ozawa has utilized a complexity function and the Dehornoy order on braids to find a canonical representation for plat presentations and bridge presentations for links \cite{ozawa2026dehornoytypeorderingplatpresentation}. See also \cite{dehornoy2008ordering} or \cite{dehornoy1994braid} for more information about the Dehornoy order on braid groups.

\medskip

Bridge presentations were introduced by Schubert in \cite{Schubert1954berEN} and have played a central role in the study of knot invariants such as bridge number. On the other hand, plat presentations were developed in braid-theoretic approaches to knot theory, particularly through the work of Birman \cite{BirmanStableEquivofPlats}. The Hilden subgroup, introduced in \cite{HildenGroup}, appears naturally as the subgroup of braid motions supported near the bridge arcs of a plat presentation.

\medskip

While the relationship between these viewpoints is often used implicitly, the equivalence between bridge isotopy and Hilden double coset equivalence does not appear in the literature in an explicit form. One contribution of this paper is to formulate and prove this correspondence in a precise way. This correspondence reduces the classification of bridge presentations up to bridge isotopy to an algebraic problem in the braid group, namely the study of Hilden double cosets.

\medskip

The paper is organized as follows. Section~\ref{sec:bridgeposition} reviews bridge presentations and bridge isotopy. Section~\ref{sec:platpresentationoflink} introduces plat presentations and the Hilden subgroup. Section~\ref{sec:TheEquivalenceofPositions} proves the main theorem. Finally, Section~\ref{sec:applicationsresults} discusses applications, including uniqueness results for special classes of knots.

\section{Bridge Position of a Link}\label{sec:bridgeposition}

Knot theory begins with the study of embeddings of $S^1$ into $S^3$ (or $\mathbb{R}^3$).  
To study knots in bridge position, we require these embeddings to satisfy a condition with respect to a height function on $S^3$.  

In general, any Morse function on $S^3$ with exactly two critical points suffices.  
However, it is standard to work with a specific choice, called the \emph{standard height function}.  
This function arises by viewing $S^3$ as the unit sphere in $\mathbb{R}^4$ and considering the projection
\[
h : \mathbb{R}^4 \to \mathbb{R}, \quad (x,y,z,w) \mapsto w.
\]
The restriction of $h$ to $S^3$ has exactly two critical points, which we denote by $\infty$ and $-\infty$.  
For each regular value $c \in \mathbb{R}$, the preimage $h^{-1}(c)$ is a $2$-sphere, called a \emph{level sphere} (or \emph{$c$-level sphere}).

Alternatively, when working in $\mathbb{R}^3$, the standard height function is given by projection onto the $z$-axis.  
In this setting, each level set $\{z = c\}$ is a plane.

\begin{figure}[ht]
    \centering
    \includegraphics[scale=0.5]{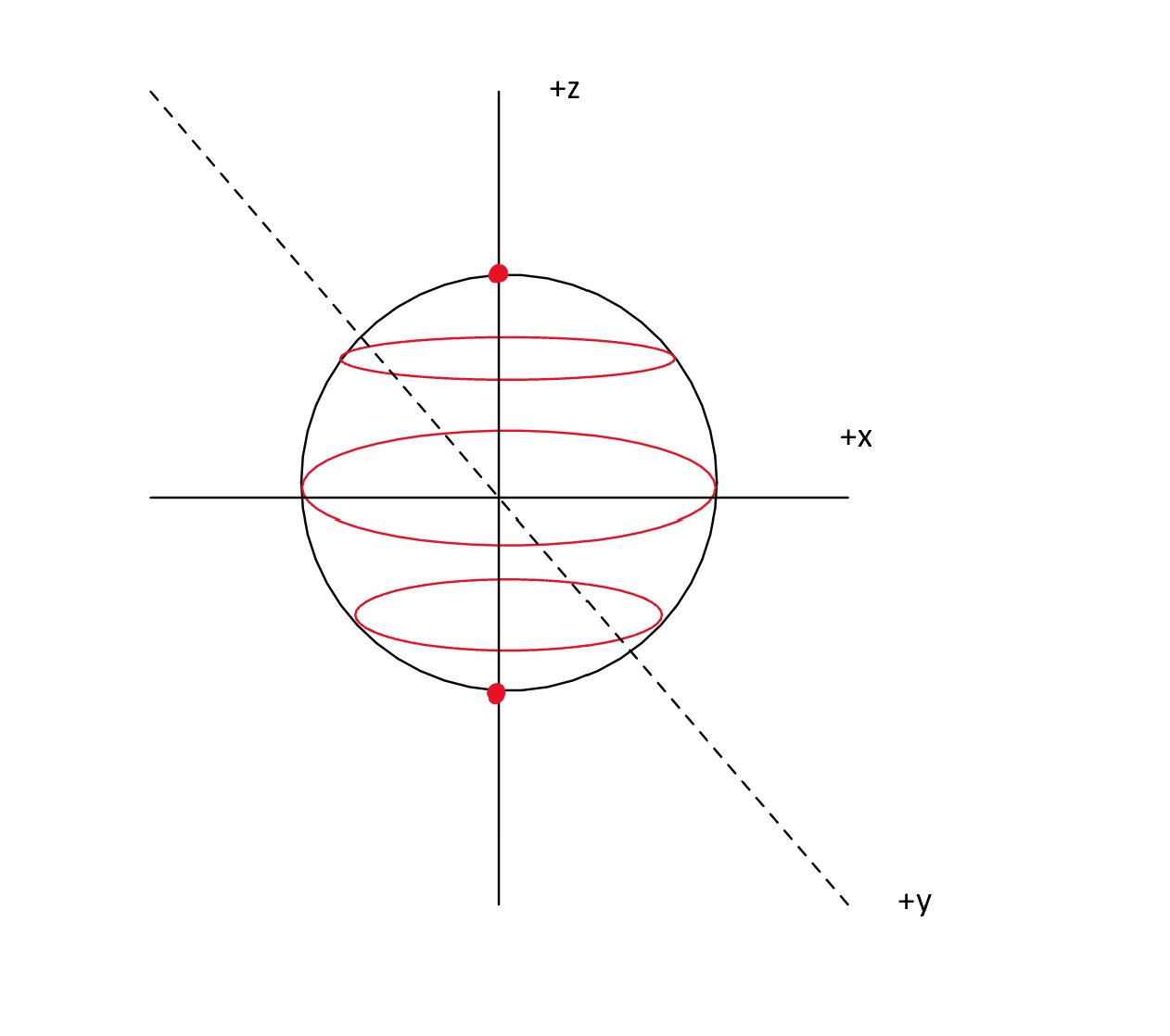}
    \caption{The height function on $S^2$ induced by projection onto the $z$-axis. Preimages of regular values are circles.}
    \label{fig:heightfunctionexample}
\end{figure}

To illustrate this idea in a lower-dimensional setting, see Figure~\ref{fig:heightfunctionexample}.

\begin{definition}
A link $L$ is in \emph{bridge position} if every local maximum of $h|_L$ lies above every local minimum.
\end{definition}

Since $L$ is a disjoint union of circles, it has the same number of local maxima and minima.  
In bridge position, these critical points are separated by level spheres.

\begin{definition}
A level $2$-sphere $S$ that separates the local minima from the local maxima is called a \emph{bridge sphere} for $L$.
\end{definition}

It is natural to consider the number of local maxima (equivalently, minima) of $h|_L$.

\begin{definition}
The \emph{bridge index} (or simply \emph{index}) of a link in bridge position is the number of local maxima.  
Equivalently, the bridge index is $n$ if the link intersects every bridge sphere in $2n$ points.  
The \emph{bridge number} of $L$, denoted $b(L)$, is the minimal bridge index over all bridge presentations of $L$.
\end{definition}

The bridge index of a specific presentation is not a link invariant, since one can introduce additional pairs of critical points. 
We now define the appropriate notion of equivalence for links in bridge position.

\begin{definition}
Two bridge presentations of a link are \emph{bridge isotopic} if there exists an ambient isotopy of $S^3$ between them such that the link remains in bridge position throughout the isotopy.
\end{definition}

By definition, bridge isotopy preserves the bridge index.  
In particular, such an isotopy cannot move a local minimum past a local maximum.  
Isotopies that violate this condition may cancel or create pairs of critical points and thus change the bridge index.

\begin{figure}[H]
    \centering
    \includegraphics[scale=0.5]{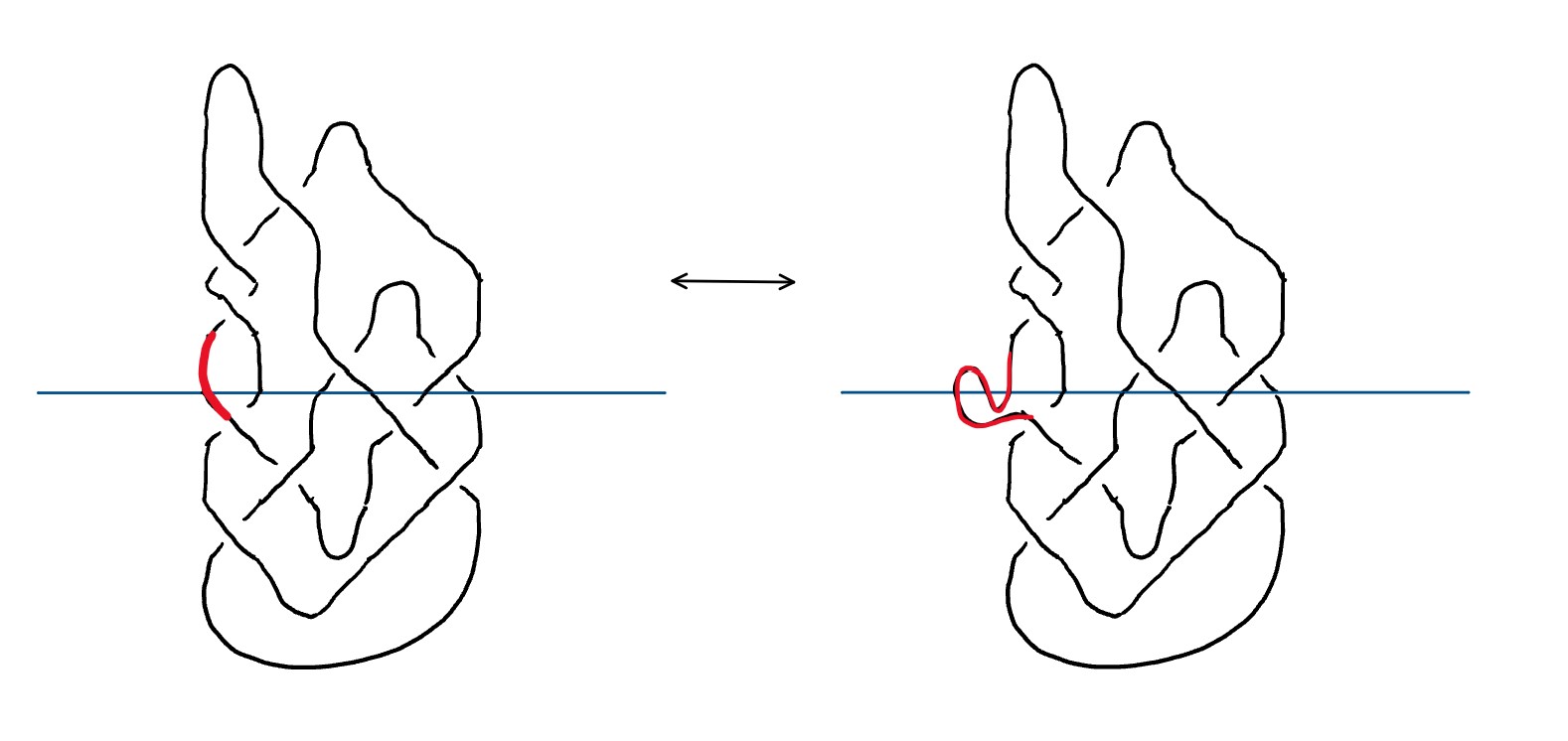}
    \caption{A link in bridge position undergoing stabilization or destabilization along the indicated subarc.}
    \label{fig:bridgestabilization}
\end{figure}

A move that cancels a local maximum and a local minimum is called a \emph{destabilization}.  
Conversely, one can introduce a pair of critical points by locally perturbing the link near a bridge sphere; this operation is called a \emph{stabilization}.  
The resulting link remains in bridge position, but its bridge index increases.  
See Figure~\ref{fig:bridgestabilization}.

The following theorem describes how different bridge presentations of the same link are related.

\begin{theorem}[\cite{Stableequiv_hayashi}]
Two links are ambient isotopic if and only if any two of their bridge presentations are related by a sequence of stabilizations, destabilizations, and bridge isotopies.
\end{theorem}

This result can be understood by viewing a bridge presentation as a knot diagram and expressing Reidemeister moves \cite{Reidemeister} as combinations of stabilizations, destabilizations, and bridge isotopies (see \cite{Stableequiv_hayashi}).

It is natural to ask whether two bridge presentations of the same link with the same bridge index must be bridge isotopic.  
In other words, one might wonder whether the bridge index is a complete invariant of bridge isotopy. This is not the case. In fact, even minimal bridge presentations (i.e., those realizing $b(L)$) need not be bridge isotopic.  
Examples of this phenomenon are discussed in Section~\ref{sec:applicationsresults}; see also \cite{montesinos_1976}.

\subsection{Bridge Positions and the Suspension of $S^2$} 
\label{subsec:Bridge Positions and the Suspension of $S^2$}

In this section, we will utilize a slightly different way of considering bridge positions, plats and their correspondence.  The main idea is we can view $S^3$ with the standard height function, as a height function coming naturally from viewing $S^3$ as the suspension of $S^2$.  
\\

Recall that the \emph{suspension} of a topological space $X,$ denoted $S X,$ is the quotient space of $X\times I$ under the identifications $(x,0)\sim (x',0)$ and $(x,1) \sim (x',1).$  These identifications collapse the subspaces $X\times \{0\}$ and $X\times \{1\}$ to a single point, we call these \emph{poles}.  Each point in $SX$ is then given by $x\in X$ and $t\in [0,1].$\\

\begin{figure}[H]
    \centering
    \includegraphics[scale=0.6]{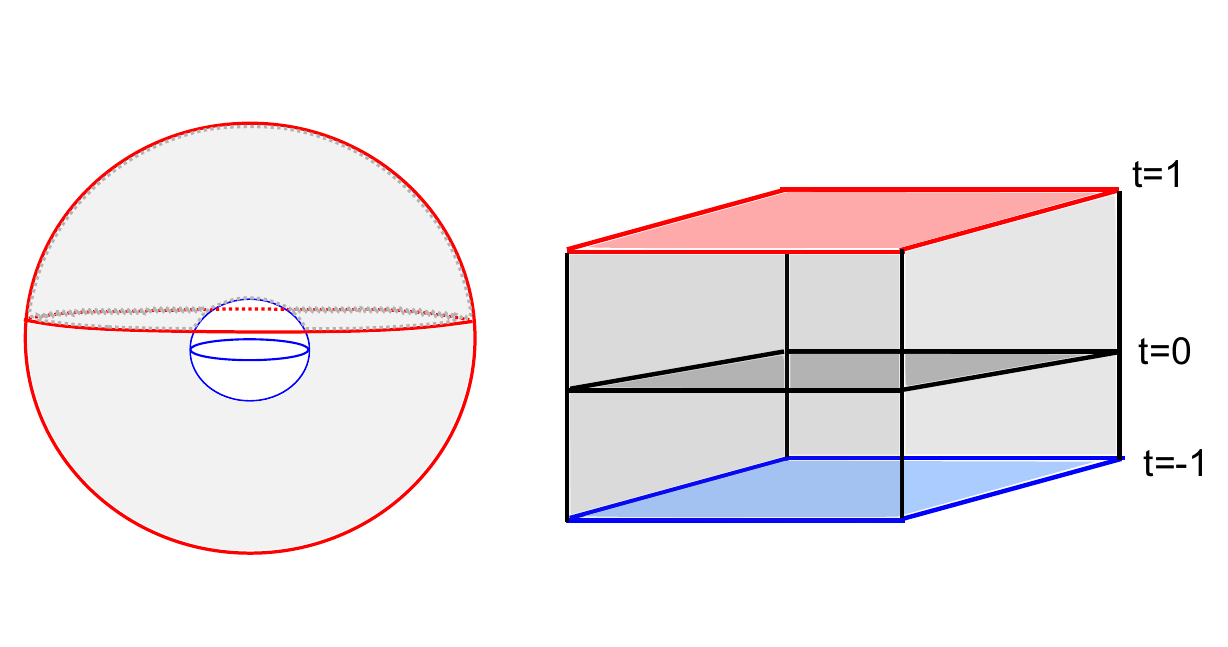}
    \caption{Suspension of $S^2.$  Left: The outermost (red) sphere and innermost (blue) sphere are contracted to a single point.  Right: Removing a ``core'' from the picture on the left we obtain a solid cube whose top and bottom face are contracted to a single point.  }
    \label{fig:suspensionofsphere}
\end{figure}

We will be considering the suspension of $S^2,$ which is homeomorphic to $S^3.$  The natural height function on $S^3$ is recovered by projecting a point in $SS^2$ to its $t$-value.  In Figure \ref{fig:suspensionofsphere}, we show $S^2 \times I$ and consider the innermost and outermost spheres to be single points.  Given a link $L$ in bridge position, we will choose a splitting sphere for $L$ and then view $S^3$ as the suspension of the splitting sphere.  Viewed this way, we will consider $S^3:= (S^2 \times [-1,1])\backslash \sim$ so that the splitting sphere is $S^2 \times \{0\}.$  All the minima of $L$ then occur on spheres with negative $t$-values and the maxima occur on spheres with positive $t$-values.

\begin{figure}[H]
    \centering
    \includegraphics[scale=0.6]{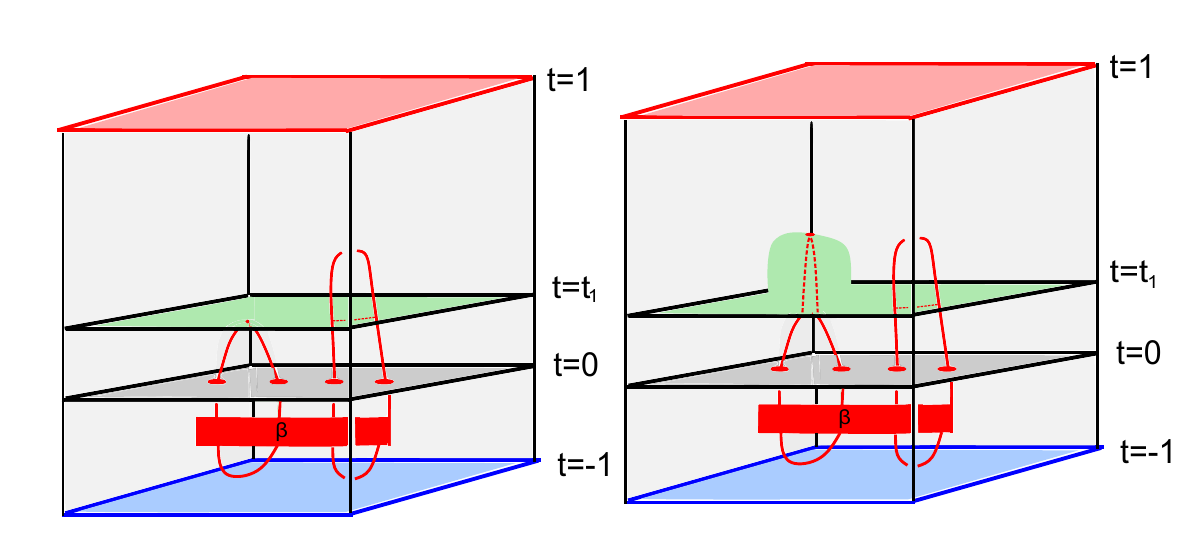}
    \caption{A bridge isotopy realized as isotoping level spheres.  The first figure is of a knot in bridge position and the $t=t_1$ sphere outlined in green. This sphere is then isotoped to push up the point where the lower bridge touches tangent to the $t=t_1$ sphere, the result on the knot is that the lower bridge is now raised, however intersections of the knot on the level spheres remains unchanged!}
    \label{fig:changeinfoliationsuspension}
\end{figure}

This perspective has several advantages. First, it provides a concrete and symmetric way to view $S^3$ as a compact manifold equipped with a height function. Second, it allows us to reinterpret a bridge isotopy not as moving the link, but as isotoping the level spheres (away from the poles) while keeping the link fixed. This viewpoint will be essential when we prove the equivalence between bridge isotopy and Hilden double coset moves. See Figure \ref{fig:changeinfoliationsuspension}.

\begin{figure}[H]
    \centering
    \includegraphics[scale=0.6]{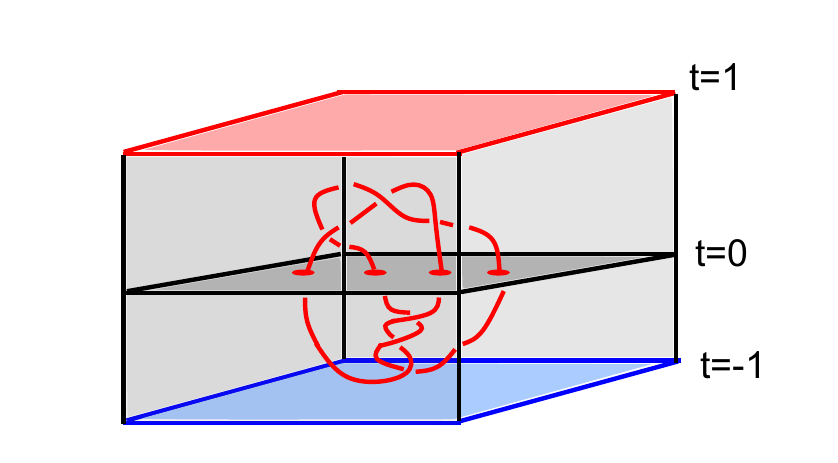}
    \caption{A knot in bridge position with respect to height function from projecting onto the $I$ value in the suspension of $S^2.$ }
    \label{fig:knotinbridgepositionsuspension}
\end{figure}

\section{Plat Presentation of a Link}\label{sec:platpresentationoflink}

A simple way to define a plat presentation is to begin with a bridge presentation and isotope the link so that all local minima occur at the same height and all local maxima occur at another.  
While this provides a valid definition, it does not explain the algebraic significance of plat presentations.  
In particular, plat closures arise naturally in braid theory; see \cite{BirmanBraidsLinksMCG, braidssurvey} for background.

\subsection{Plats and the Hilden Double Coset}

Given a braid word $\beta \in B_{2n}$, a \emph{plat presentation} is obtained by connecting the $i$-th strand to the $(i+1)$-th strand at both the top and bottom, for each odd index $i$.  We provide an example below.  

\begin{example}
\text{}\\
\begin{figure}[h]
    \subcaptionbox{The braid $\beta = \sigma_2 \sigma_4 \sigma_1 \sigma_3 \sigma_1$.}[0.45\textwidth]{
    \centering
    \begin{tikzpicture}[scale=0.6]
    \pic[scale=0.6,braid/.cd,
        number of strands=6,
        line width=2pt,
        name prefix=braid] at (2,0) {braid={s_2 s_4 s_1 s_3 s_1}};
    \end{tikzpicture}
    }
    \hfill
    \subcaptionbox{The plat closure of $\beta$}[0.45\textwidth]{
    \centering
    \hspace*{2.2cm}
    \begin{tikzpicture}[scale=0.6]
        \pic[scale=0.6,braid/.cd,
         number of strands=6,
        line width=2pt,
        name prefix=braid, 
        rotate=0,
        ] at (2,0) {braid={| s_2 s_4 s_1 s_3  s_1}};
        \node[fill=red,circle,inner sep=0pt,minimum size=0pt] at (2, 0)  (t1) {};
        \node[fill=red,circle,inner sep=0pt,minimum size=0pt] at (3, 0)  (t2){};
        \node[fill=red,circle,inner sep=0pt,minimum size=0pt] at (4, 0)  (t3){};
        \node[fill=red,circle,inner sep=0pt,minimum size=0pt] at (5, 0)  (t4){};
        \node[fill=red,circle,inner sep=0pt,minimum size=0pt] at (6, 0)  (t5){};

        \node[fill=red,circle,inner sep=0pt,minimum size=0pt] at (7, 0)  (t6){};
        \node[fill=red,circle,inner sep=0pt,minimum size=0pt] at (2, -5.5)  (b1){};
        \node[fill=red,circle,inner sep=0pt,minimum size=0pt] at (3, -5.5)  (b2){};
        \node[fill=red,circle,inner sep=0pt,minimum size=0pt] at (4, -5.5)  (b3){};
        \node[fill=red,circle,inner sep=0pt,minimum size=0pt] at (5, -5.5)  (b4){};
        \node[fill=red,circle,inner sep=0pt,minimum size=0pt] at (6, -5.5)  (b5){};

        \node[fill=red,circle,inner sep=0pt,minimum size=0pt] at (7, -5.5)  (b6){};
    
        \draw[black,line width=2pt, bend left]  (t1) to node [auto] {} (t2);
         \draw[black,line width=2pt, bend left]  (t3) to node [auto] {} (t4);
         \draw[black,line width=2pt, bend left]  (t5) to node [auto] {} (t6);
         \draw[black,line width=2pt, bend right]  (b1) to node [auto] {} (b2);
         \draw[black,line width=2pt, bend right]  (b3) to node [auto] {} (b4);
          \draw[black,line width=2pt, bend right]  (b5) to node [auto] {} (b6);
        \end{tikzpicture}
    }
\caption{Plat closures of braids in $B_{2n}$}
\label{fig:platclosure}
\end{figure}
\end{example}

The plat above in Figure \ref{fig:platclosure} is called a \emph{$6$-plat}, since $\beta \in B_6$, or equivalently a \emph{$3$-bridge plat}.  
Notice that any plat presentation is already in bridge position.

\medskip

A fundamental problem is to determine the link type represented by a given plat.  
In general, this is quite difficult: complicated braid words may yield trivial links (see Figure~\ref{fig:hardunknots}).  

A more tractable question is the following: given a plat presentation of a fixed link type, which operations on the braid word preserve that link type?  
One approach is to interpret the plat as a knot diagram and apply Reidemeister moves that preserve the plat structure.  
These moves can then be translated into operations in the braid group.

We now describe these moves.

\subsubsection{Reidemeister Moves on Plat Diagrams}

\begin{figure}[H]
\centering
\includegraphics[scale=0.7]{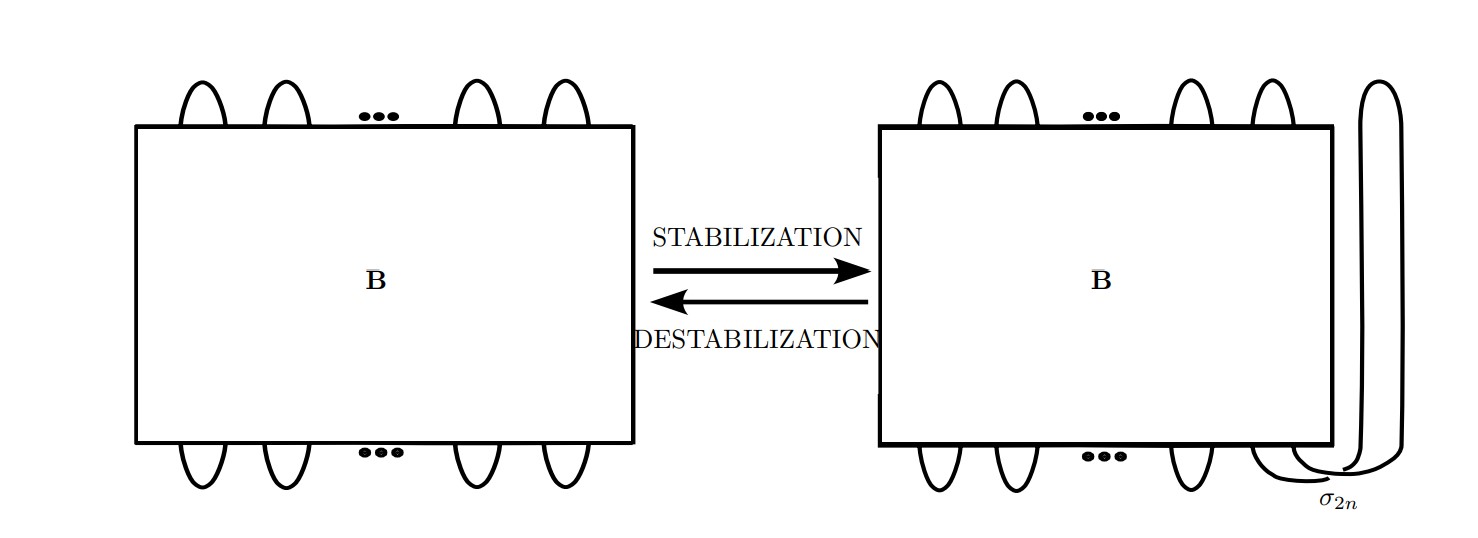}
\caption{Stabilization and destabilization of a plat.}
\label{fig:Platstabilization}
\end{figure}

\begin{figure}[H]
\centering
\includegraphics{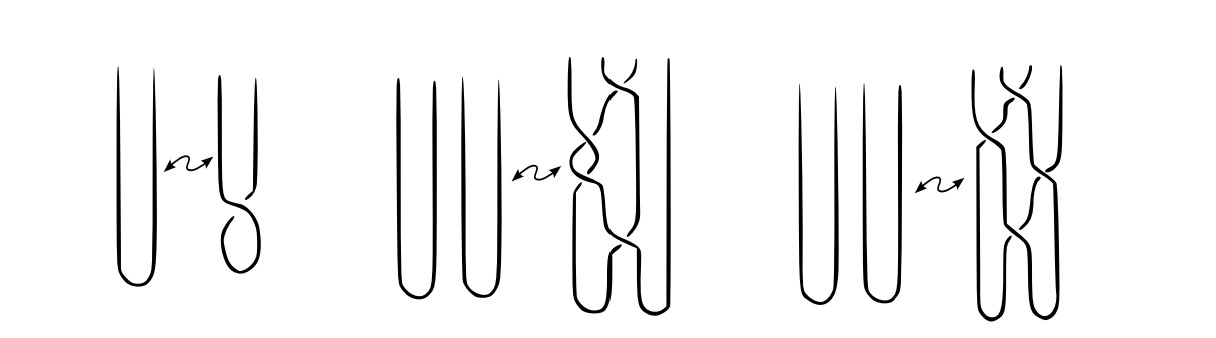}
\caption{Examples of Type I and Type II moves on a plat diagram.}
\label{fig:reidmoves}
\end{figure}

\begin{enumerate}
    \item \textbf{Type I moves.}  
    To preserve the plat structure, Type I moves can only be applied in neighborhoods containing a bridge.  
    Performing such a move elsewhere would introduce a new critical point at a different height, violating the definition of a plat.

    \item \textbf{Type II moves.}  
    These can be applied away from the bridges without introducing new critical points.  
    When applied near bridges, they allow one to move bridges past one another (see Figure~\ref{fig:reidmoves})\footnote{Figures \ref{fig:reidmoves} and \ref{fig:Platstabilization} are borrowed with permission from \cite{solanki2023studyinglinksplatsunlink}}.

    \item \textbf{Type III moves.}  
    These correspond exactly to braid isotopies.  
    They change the order of crossings without affecting the bridge structure.

    \item \textbf{Stabilization.}  
    A Type I move performed away from the bridges creates a new pair of critical points.  After isotoping back into plat position, this yields a \emph{stabilization}.  
    This operation increases the number of strands.  
    Birman \cite{BirmanStableEquivofPlats} showed that stabilization on any strand is equivalent to stabilization on an outermost strand, up to braid isotopy and bridge moves.
\end{enumerate}

The moves described above that occur near the bridges generate a subgroup of the braid group known as the \emph{Hilden subgroup}, denoted $H_{2n}$.  We will call these \emph{Hilden double coset moves.}  Alternatively,
let $\tau_n\subset B^3$ denote the standard trivial $n$–string tangle, given by the
$n$ disjoint, unknotted bridges of a plat presentation. Then the Hilden subgroup may be defined as the braids whose induced action on the $2n$ boundary points extends to a homeomorphism of the pair 
$(B^3,\tau_n).$  This group was also studied in \cite{wicketsbrendlehatcher}, where it appears as the space of disjoint wickets.

\begin{definition}[Hilden subgroup]\label{Defn:Hildensubgroup}
An element $\gamma \in B_{2n}$ lies in the \emph{Hilden subgroup} $H_{2n}$ if, for every $\beta \in B_{2n}$, the plat closures of $\gamma\beta$ and $\beta\gamma$ represent the same link type as the plat closure of $\beta$.  
The group $H_{2n}$ is generated by
\[
\{\sigma_1,\ \sigma_2\sigma_1^2\sigma_2,\ \sigma_{2i}\sigma_{2i-1}\sigma_{2i+1}\sigma_{2i} \mid 1 \leq i \leq n-1\}.
\]
\end{definition}

These generators correspond precisely to the bridge moves described above (see Figure~\ref{fig:reidmoves}).  
A presentation for $H_{2n}$ was given by Tawn \cite{tawnpresentation}.

\medskip

Multiplying a braid word on the left or right by elements of $H_{2n}$ does not change the link type of its plat closure.  
Thus, braid words in the same \emph{Hilden double coset} determine the same link type.

This leads to two complementary viewpoints:
\begin{itemize}
    \item \emph{Algebraic:} two braid words are equivalent if they lie in the same Hilden double coset.
    \item \emph{Geometric:} two plat diagrams are equivalent if they are related by bridge moves and isotopy.
\end{itemize}

\begin{definition}
Two $2n$-plats are \emph{plat equivalent} if they are related by Type I - III moves described above.  
Equivalently, they lie in the same Hilden double coset.
\end{definition}

\begin{figure}[H]
\centering
\includegraphics[scale=0.5]{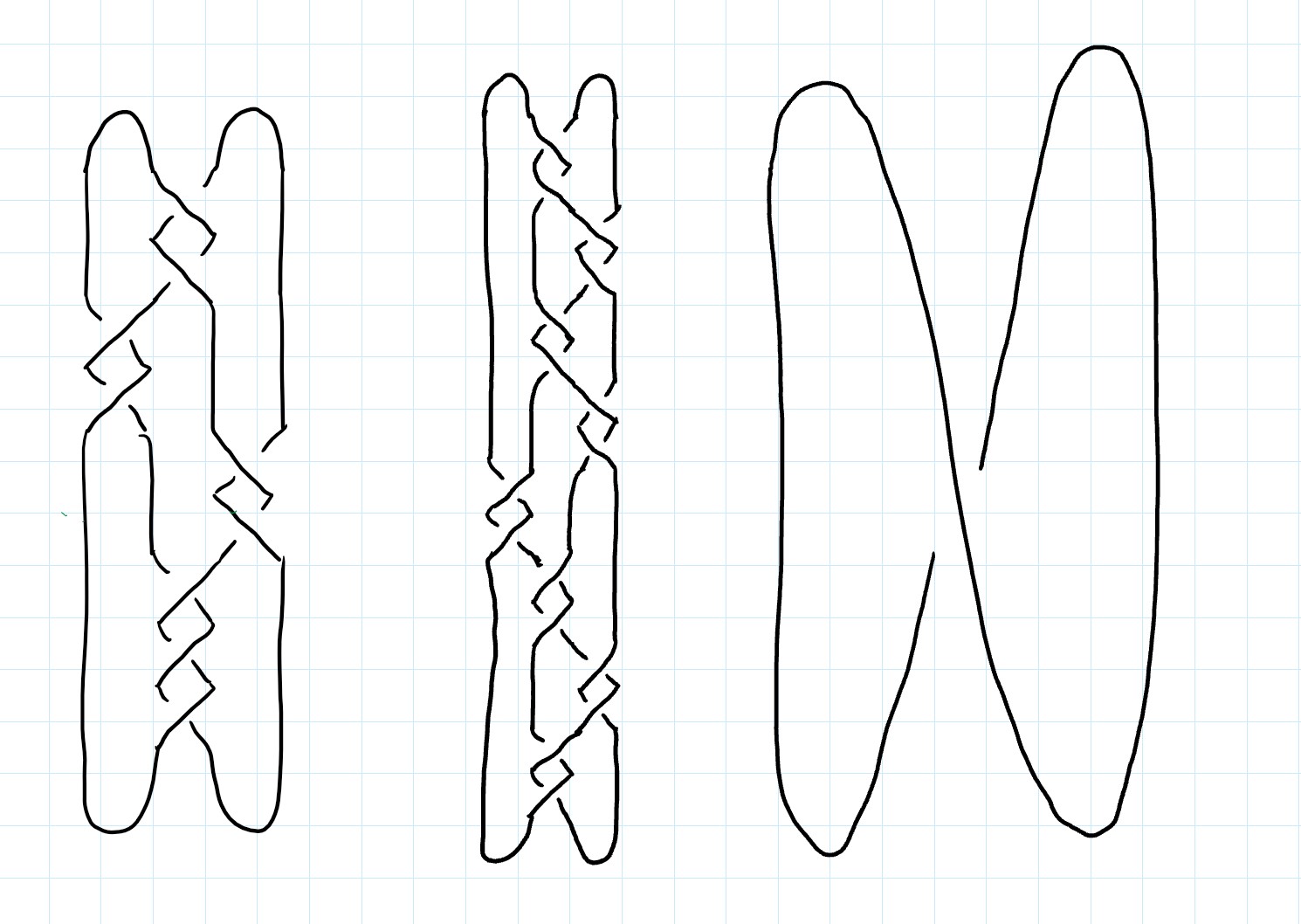}
\caption{Different plat presentations of the unknot in the same double coset.}
\label{fig:hardunknots}
\end{figure}

Examples of distinct plat diagrams representing the same link type are shown in Figure~\ref{fig:hardunknots}.

\begin{definition}
Two $2n$-plats are \emph{stably equivalent} if they are related by Type I-III moves together with stabilization and destabilization.
\end{definition}

Birman \cite{BirmanStableEquivofPlats} proved the following analogue of Markov's theorem.

\begin{theorem}[Markov's Theorem for Plats]
All plat presentations of a fixed link type are related by braid isotopies, bridge moves, and stabilizations/destabilizations.
\end{theorem}

\subsection{Alternative View of Plat Position}

Using the suspension model of $S^3$ described in Section~\ref{subsec:Bridge Positions and the Suspension of $S^2$}, we may also view plat presentations in terms of the height function $t \in [-1,1]$.  

A link is in \emph{plat position} if all local maxima occur at a single level $t_{\max}$ and all local minima occur at a single level $t_{\min}$.  
In this setting, it is useful to track how bridge moves act on the intersection points of the link with the level spheres at $t_{\max}$ and $t_{\min}$.  

\begin{figure}[H]
\centering
\includegraphics[scale=0.8]{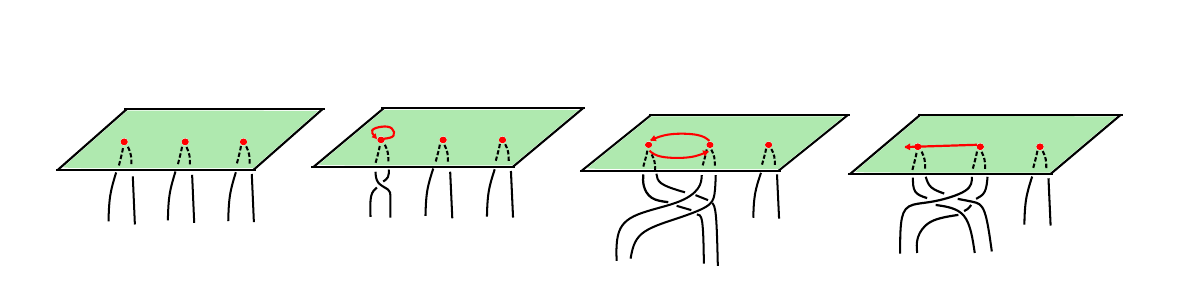}
\caption{Bridge moves viewed as homotopies on a level sphere.}
\label{fig:bridgemovesonmaxlevel}
\end{figure}

These moves can be interpreted as homotopies of points on the level spheres, as illustrated in Figure~\ref{fig:bridgemovesonmaxlevel}.  
This perspective will be important in establishing the equivalence between plat presentations and bridge positions in Section~\ref{sec:TheEquivalenceofPositions}.

\section{Equivalence of Plat Presentations and Bridge Isotopy}\label{sec:TheEquivalenceofPositions}
In this section, we show the main result of this paper.  We say that plat presentations $P$ and $P'$ \emph{correspond} to bridge presentations $L$ and $L'$ if they are obtained by isotoping all local maxima and minima of $L$ and $L'$ to common levels without changing the number of bridges.

\begin{theorem}\label{maintheorem}
 If two bridge presentations $L$ and $L'$ are bridge isotopic if and only if any plats $P$ and $P'$ corresponding to $L$ and $L'$ respectively, are in the same Hilden double coset class.  Thus, the set of bridge isotopy classes of
$n$-bridge presentations is in bijection with Hilden double cosets in $B_{2n}.$
\end{theorem}

\begin{figure}[ht]
        \centering
        \includegraphics[width=\linewidth]{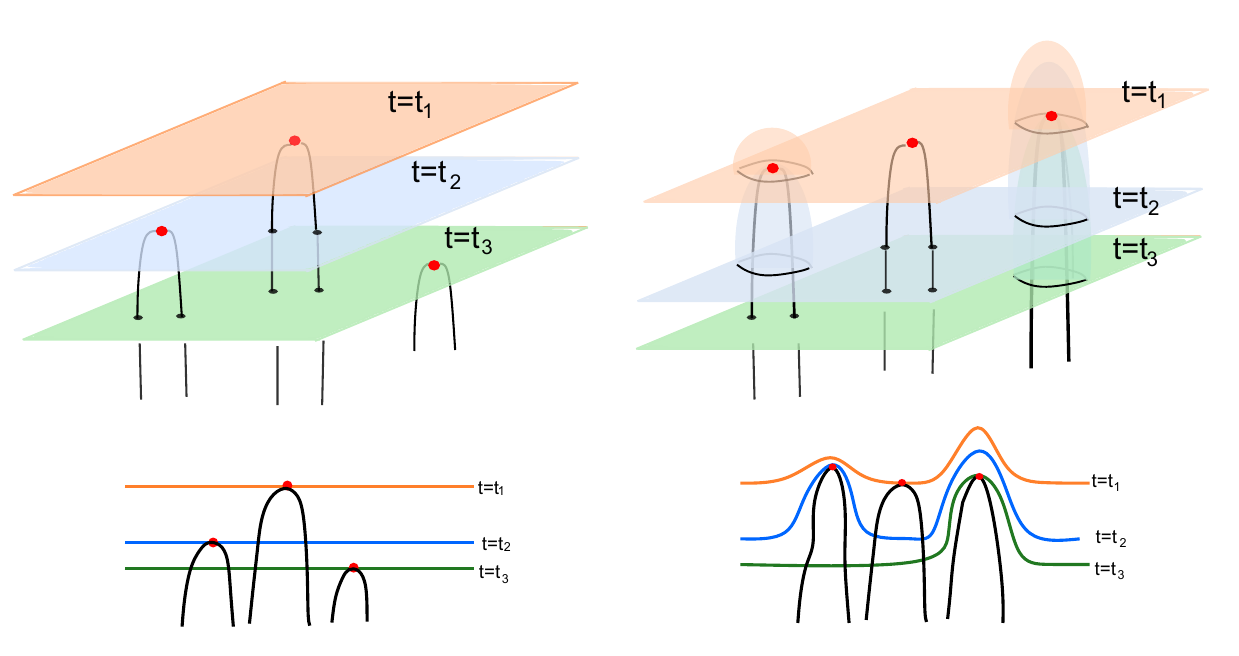}
        \caption{Pushing the level pages up so that each local maximum appears to be at the same height}
        \label{fig:changeinfoliationtoplatpostition}
    \end{figure}

\begin{figure}[ht]
    \centering
    \includegraphics[scale=0.7]{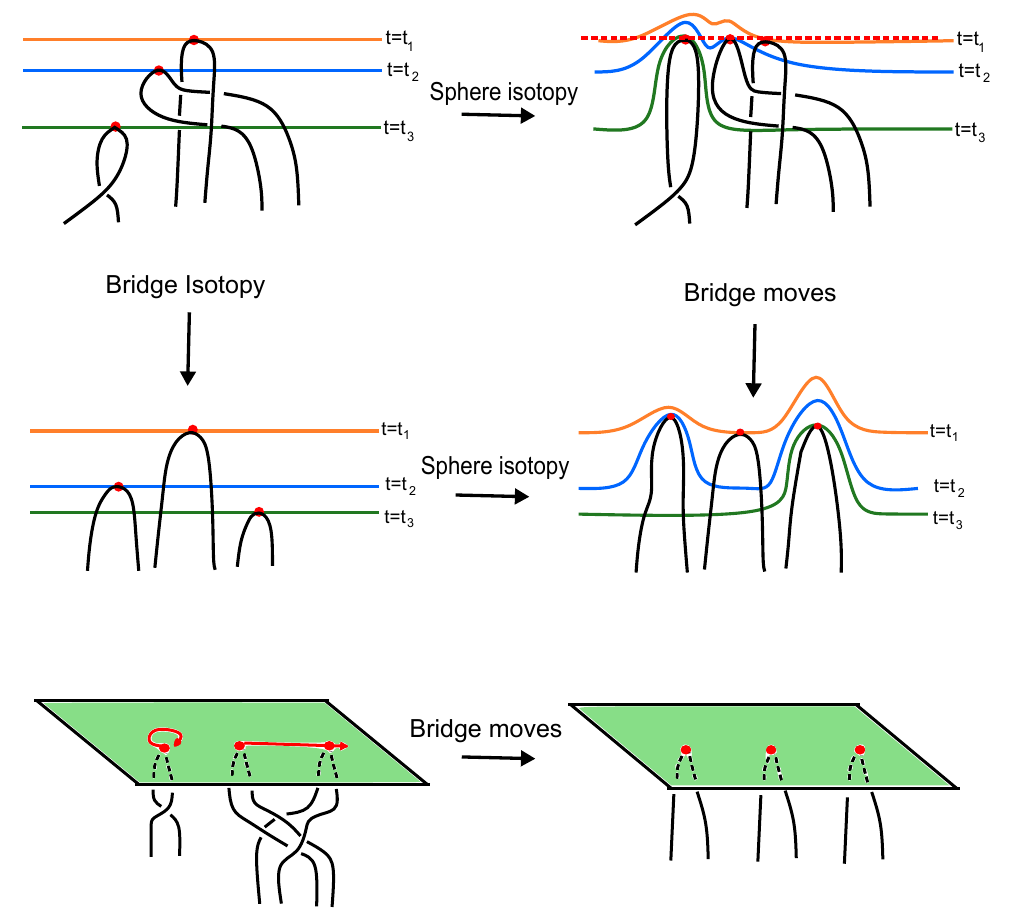}
    \caption{The top four pictures show $L$ and $L'$ and the isotopy of level spheres that take them to $P$ and $P'$ (moving from left to right).  The bottom two pictures show how the bridge isotopy of the new collection of spheres taking $L$ to $L'$ viewed from the dotted level sphere (see top right figure) is exactly a sequence of bridge moves.}
    \label{fig:lengthyexample}
\end{figure}

\begin{proof}

First suppose that $P$ and $P'$ lie in the same Hilden double coset. Then there exist $h_1, h_2 \in H_{2n}$ such that
\[
P' = h_1 \, P \, h_2.
\]
By definition of the Hilden subgroup, multiplication by $h_1$ and $h_2$ corresponds to performing bridge moves and braid isotopies. These moves preserve the plat structure and are realized by bridge isotopies. Therefore, $P$ and $P'$ represent bridge isotopic links.

\medskip

For the converse, suppose that $L$ and $L'$ are bridge isotopic. We first show that, after isotoping the level spheres, both links can be viewed in plat position.

View $S^3$ as the suspension of a splitting sphere, with height function given by projection onto the $t$-coordinate. Suppose the local maxima of $L$ occur at levels $t = t_i$, for $i = 1,\dots,n$. Each level sphere $t = t_i$ contains a point of tangency with the link.

Starting from the lowest maximum, isotope a small neighborhood of each tangency upward so that all maxima occur at the same level. This isotopy is supported in disjoint neighborhoods and avoids the rest of the link. Repeating this process for all maxima, and similarly for minima, we obtain a new foliation of $S^3$ in which $L$ appears in plat position (see Figure~\ref{fig:changeinfoliationtoplatpostition}).

Importantly, this process does not change how $L$ intersects each level sphere; it only modifies the foliation. Thus, the link itself remains unchanged, but is now viewed as a plat $P$ with respect to the new height function.

Applying the same construction to $L'$ produces a plat $P'$. 
Although there is some ambiguity in the choice of isotopies of the critical points to the maximum and minimum the level spheres, any two such resulting foliations are isotopic through foliations fixing the poles, and hence correspond to multiplication by elements of the Hilden subgroup (see Figures \ref{fig:bridgemovesonmaxlevel} and \ref{fig:lengthyexample}).

Since $L$ and $L'$ are bridge isotopic, this isotopy can be realized with respect to the modified foliations. Viewing the isotopy in this way, it induces a sequence of moves that preserve the bridge structure at the top and bottom, and hence correspond precisely to multiplication by elements of the Hilden subgroup. Therefore, $P$ and $P'$ lie in the same Hilden double coset.

\end{proof}

  Figure~\ref{fig:bridgemovesonmaxlevel}.  Figure \ref{fig:lengthyexample} below illustrates how the bridge isotopy looks and the corresponding plat isotopies. 

\begin{corollary}[Markov's Theorem for Plats]
   If $P$ and $P'$ are equivalent as links, they are related by a series of braid isotopies, (de)stabilizations, and bridge moves.
\end{corollary}
\begin{proof}
    All bridge positions of a fixed link type are related by bridge isotopy and (de) stabilizations.  By Theorem \ref{maintheorem}, $P$ and $P'$ are related by (de) stabilization and bridge isotopy.  Therefore, after stabilizing or destabilizing a suitable number of times they are bridge isotopic, or in the same Hilden double coset class.
\end{proof}

\section{Applications of Equivalence}\label{sec:applicationsresults}
As stated above, the importance of the main theorem is that it ties together two ways of considering links in $S^3.$  We outline this with a few examples. 
\subsection{The Hilden Double Coset Problem}
It was shown by Alexander that any link type can be represented as a closed braid.  Markov then showed that if two braid closures represent the same link type then they are related by braid isotopy, (de) stabilization, and conjugation in the braid group.  For more information on representing links as braid closures see \cite{braidssurvey}.  The powerful advantage of representing a link as a braid closure, as opposed to a plat closure, is that both the word and conjugacy problems have been solved in the braid groups.  That is, one can algorithmically determine whether two braid words are isotopic or if they are conjugate to each other.  The downside to this way of representing links is that there are braid closures of words in $B_3$ that represent the unknot that are not conjugate to each other.  In fact, there are infinitely many conjugacy classes of the unknot. This is a good motivation for studying plat closures! It is worthwhile to investigate another equivalence relation on braids where the equivalence classes are ``larger'' than conjugacy classes.  In fact we will show below, that there is one Hilden double coset class of the $n$-bridge unknot. The dilemma we are faced with is that there is currently no algorithmic way to determine whether two braid words are in the same Hilden double coset class.  Indeed, algebraically considered this is a very hard problem \cite{kalka2015doublecosetproblemparabolic}. 
 This brings us to the first application of plat equivalence and bridge equivalence:

\begin{theorem}
    Two braid words $\alpha,\beta\in B_{2n}$ are in the same Hilden double coset class if and only if their plat closures are bridge isotopic.
\end{theorem}

This means that a method for determining if two links in bridge position are bridge isotopic would provide a solution to the Hilden double coset problem.

\subsection{Summary of Results about Hilden Double Coset Classes of Special Knots}
Even with the double coset problem as a stumbling block, plat closures do provide some very interesting results.  The problem of having infinitely many conjugacy classes of braids that have the same braid closure does not arise so readily when taking the plat closure. We consider Hilden double cosets of special knots in this section. To start, we consider the unknot.  
It has been shown by Otal, Schultens, and others \cite{Schultens_2009}\cite{Otal} that: 
\begin{theorem}
    All bridge presentations of the unknot of a given bridge index are equivalent.  
\end{theorem}

Our correspondence shows then that: 
\begin{theorem}
    There is only one Hilden double coset class of the unknot.
\end{theorem}

This was also proved by Solanki without the correspondence above in \cite{solanki2023studyinglinksplatsunlink}.
It was shown in \cite{TorusknotsoneHDC} there is a single bridge isotopy class of any $n$-bridge torus knot as well.  This gives:
\begin{theorem}
    There is only one Hilden double coset class of torus knots
\end{theorem}

Additionally, it was shown in \cite{BirmanBraidsLinksMCG} there are at most four distinct double coset classes for two bridge plats.  This was curious to the author since we could only find two distinct classes.  The reason for this is answered by the correspondence: 
\begin{theorem}\cite{twobridgeflip}
    There are at most two distinct isotopy classes for two bridge knots.  These are related by a map that reflects the knot and takes the top bridges to the bottom.
\end{theorem}

 Menasco and Solanki \cite{menasco2024studyinglinksplatssplit} considered links in plat position and showed that:
\begin{theorem}(\cite{menasco2024studyinglinksplatssplit})
    In each Hilden double coset class of a split link, there is a representative that is obviously split. 
\end{theorem}

This means that in every bridge isotopy class of a split link there is a representative that is obviously split. 
\\
In \cite{ozawa_takao_2013} it was shown that there exists a 4-bridge position of a 3-bridge knot that does not admit a destabilization. Pulling this bridge position to a plat position then gives a 4-bridge plat presentation of a 3-bridge knot that does not admit a destabilization.  We summarize this result below as a theorem: 
\begin{theorem}
There exist plat presentations whose double coset class does not admit a destabilization but the bridge index is strictly less than the given plat presentations index. 
\end{theorem}

In \cite{Schultens_2009} Schultens asks: Are there prime knots with width complexes containing distinct global minima?  The question framed in terms of plat presentations is: Are there distinct Hilden double coset classes of prime knots that realized the bridge number?  The answer is yes.  Montesinos in \cite{montesinos_1976} gave an example of a knot that had two distinct minimal plat presentations, by our equivalence we have two distinct minimal bridge presentations.

\section{Conclusion and Further Discussion}

The results of this paper serve as a bridge between two ways of viewing links.  The geometric power of viewing a link in bridge position is supplemented by the algebraic approach of representing a link as a plat.  It is no doubt that the current literature on these two presentations favors considering bridge positions of links. In forthcoming work, the author will explore developing invariants and combinatorial objects that arise from representing links as plats.  By the results of this current paper, these invariants and combinatorial objects will naturally apply to links in bridge position.\\ 
\text{}\\
 We end this section by outlining a couple of interesting directions for further research.
\begin{question}
    For a 3-bridge knot are there only finitely many non-isotopic minimal bridge positions? 
\end{question}

\begin{question}
    Can you find two 3-bridge presentations of a fixed knot type that are not bridge isotopic and require more than two stabilizations to become isotopic? This is analogous to the Heegaard splitting stabilization conjecture in genus 2.  
\end{question}

\newpage

\bibliography{updatedbib} 
\vspace{15mm}

\newcommand{\Addresses}{{
  \bigskip
    Seth Hovland, \textsc{Department of Mathematics, University of Minnesota Duluth, USA}\par\nopagebreak
  \textit{E-mail address}: \texttt{hovl0070@d.umn.edu}

}}

\Addresses

\end{document}